\documentclass[english]{amsart}
\usepackage[T1]{fontenc}
\usepackage[latin9]{inputenc}
\usepackage{mathrsfs}
\usepackage{amsthm}
\usepackage{amstext}
\usepackage{amssymb}
\usepackage{esint}

\makeatletter

\newcommand{\noun}[1]{\textsc{#1}}

\theoremstyle{plain}
\newtheorem{thm}{\protect\theoremname}
  \theoremstyle{plain}
  \newtheorem{cor}[thm]{\protect\corollaryname}

\makeatother

\usepackage{babel}
  \providecommand{\corollaryname}{Corollary}
\providecommand{\theoremname}{Theorem}

\begin{document}
\global\long\def\Zp{\mathbb{Z}_{p}}

\title{Identities of symmetry for higher-order $q$-Euler polynomials}

\author{Dae San Kim \and Taekyun Kim}
\begin{abstract}
In this paper, we derive basic identities of symmetry in two variables
related to higher-order $q$-Euler polynomials and $q$-analogue of
higher-order alternating power sums. The derivation of identities
are based on the multivariate $p$-adic fermionic integral expression
of the generating function for the higher-order $q$-Euler polynomials.
\end{abstract}
\maketitle

\section{Introduction}

$\,$

Let $p$ be a fixed prime number. Throughout this paper, $\Zp$, $\mathbb{Q}_{p}$
and $\mathbb{C}_{p}$ will, respectively, denote the ring of $p$-adic
rational integers, the field of $p$-adic rational numbers and the
completion of algebraic closure of $\mathbb{Q}_{p}$.

Let $v_{p}$ be the normalized exponential valuation of $\mathbb{C}_{p}$
with $\left|p\right|_{p}=p^{-v_{p}\left(p\right)}=\frac{1}{p}$. When
one talks about $q$-extension, $q$ is variously considered as an indeterminate,
a complex number $q\in\mathbb{C}$ or $p$-adic number $q\in\mathbb{C}_{p}$.
If $q\in\mathbb{C}$, one usually assumes $\left|q\right|<1$;
if $q\in\mathbb{C}_{p}$, one usually assumes $\left|1-q\right|_{p}<1$.
The $q$-number of $x$ is defined by $\left[x\right]_{q}=\frac{1-q^{x}}{1-q}$.
Note that ${\displaystyle \lim_{q\rightarrow1}\left[x\right]_{q}=x}.$

For $r\in\mathbb{N}$, the Euler polynomials of order $r$ are defined
by the generating function to be
\begin{equation}
\left(\frac{2}{e^{t}+1}\right)^{r}e^{xt}=\left(\frac{2}{e^{t}+1}\right)\times\cdots\times\left(\frac{2}{e^{t}+1}\right)e^{xt}=\sum_{n=0}^{\infty}E_{n}^{\left(r\right)}\left(x\right)\frac{t^{n}}{n!}.\label{eq:1}
\end{equation}

When $x=0$, $E_{n}=E_{n}\left(0\right)$ are called the Euler numbers
of order $r$ (see \cite{key-4,key-6}).

From (\ref{eq:1}), we note that
\begin{equation}
E_{n}^{\left(r\right)}\left(x\right)=\sum_{l=0}^{n}\dbinom{n}{l}E_{l}^{\left(r\right)}x^{n-l},\label{eq:2}
\end{equation}
(see \cite{key-4,key-6}).

In \cite{key-4}, Kim considered the $q$-extension of (\ref{eq:1})
which is given by

\begin{equation}
2^{r}\sum_{m=0}^{\infty}\left(-1\right)^{m}\dbinom{m+r-1}{m}e^{\left[m+x\right]_{q}t}=\sum_{n=0}^{\infty}E_{n,q}^{\left(r\right)}\left(x\right)\frac{t^{n}}{n!}.\label{eq:3}
\end{equation}

Thus, by (\ref{eq:3}), we get
\begin{equation}
\sum_{n=0}^{\infty}E_{n,q}^{\left(r\right)}\left(x\right)\frac{t^{n}}{n!}=2^{r}\sum_{m_{1},\cdots,m_{r}=0}^{\infty}\left(-1\right)^{m_{1}+\cdots+m_{r}}e^{\left[m_{1}+\cdots+m_{r}+x\right]_{q}t},\label{eq:4}
\end{equation}
where $E_{n,q}^{\left(r\right)}\left(x\right)$ are called the $q$-Euelr
polynomials of order $r$$\left(\in\mathbb{N}\right)$.

When $x=0$, $E_{n,q}^{\left(r\right)}=E_{n,q}^{\left(r\right)}\left(0\right)$
are called the $q$-Euler number of order $r$.

From (\ref{eq:4}), we note that
\begin{equation}
E_{n,q}^{\left(r\right)}\left(x\right)=\sum_{l=0}^{n}\dbinom{n}{l}q^{lx}E_{l,q}^{\left(r\right)}\left[x\right]_{q}^{n-l}.\label{eq:5}
\end{equation}

Let $\mathscr{C}\left(\Zp\right)$ be the space of continuous functions
on $\Zp$. For $f\in\mathscr{C}\left(\Zp\right)$, the fermionic $p$-adic
integral on $\Zp$ is defined by Kim as follows :

\begin{equation}
I_{-1}\left(f\right)=\int_{\Zp}f\left(x\right)d\mu_{-1}\left(x\right)=\lim_{N\rightarrow\infty}\sum_{x=0}^{p^{N}-1}f\left(x\right)\left(-1\right)^{x},\label{eq:6}
\end{equation}
(see \cite{key-1,key-7,key-12}).

By (\ref{eq:6}), we easily geet
\begin{equation}
\int_{\Zp}f\left(x+n\right)d\mu_{-1}\left(x\right)+\left(-1\right)^{n-1}\int_{\Zp}f\left(x\right)d\mu_{-1}\left(x\right)=2\sum_{l=0}^{n-1}\left(-1\right)^{n-1-l}f\left(l\right),\label{eq:7}
\end{equation}
where $n\in\mathbb{N}$ (see \cite{key-2,key-3,key-5,key-9,key-10,key-11}).

Thus, from (\ref{eq:7}), we have
\begin{equation}
\int_{\Zp}e^{\left(x+y\right)t}d\mu_{-1}\left(y\right)=\frac{2}{e^{t}+1}e^{xt}=\sum_{n=0}^{\infty}E_{n}\left(x\right)\frac{t^{n}}{n!}.\label{eq:8}
\end{equation}

By (\ref{eq:8}), we easily get
\begin{align}
\int_{\Zp}\cdots\int_{\Zp}e^{\left(x+y_{1}+\cdots+y_{r}\right)t}d\mu_{-1}\left(y_{1}\right)\cdots d\mu_{-1}\left(y_{r}\right) & =\left(\frac{2}{e^{t}+1}\right)^{r}e^{xt}\label{eq:9}\\
 & =\sum_{n=0}^{\infty}E_{n}^{\left(r\right)}\left(x\right)\frac{t^{n}}{n!}.\nonumber
\end{align}

In the next section, we consider the $q$-analogue of (\ref{eq:9}).
In \cite{key-3}, the thirty one basic identities of symmetry in three
variables related to higher-order Euler polynomials and alternating
power sums are derived from (\ref{eq:9}).

In \cite{key-5}, Kim gave some interesting relations of symmetry
between the alternating power sum polynomials and Euler polynomials
and he suggested an open question as to finding the $q$-extension of symetry
$p$-adic invariant integral on $\Zp$ for $q$-Euler polynomials.

Recently, several authors have studied the identity of symmetry and
$q$-extensions of Euler polynomials which are derived from the $p$-adic
fermionic integrals on $\Zp$ (see {[}1-13{]}).

In this paper, we investigate several further interesting properties
of symmetry for the multivariate $p$-adic fermionic integrals on $q$-polynomials.
From our investigation, we derive some relations of symmetry between
the higher-order alternating power sum $q$-polynomials and the higher-order
$q$-Euler polynomials.

\section{Identities of symmetry for the higher-order $q$-Euler polynomials}

$\,$

From (\ref{eq:3}), (\ref{eq:4}), (\ref{eq:6}) and (\ref{eq:7}),
we note that
\begin{align}
 & \int_{\Zp}\cdots\int_{\Zp}e^{\left[x+y_{1}+\cdots+y_{r}\right]_{q}t}d\mu_{-1}\left(y_{1}\right)\cdots d\mu_{-1}\left(y_{r}\right)\label{eq:10}\\
= & 2^{r}\sum_{m_{1},\cdots,m_{r}=0}^{\infty}\left(-1\right)^{m_{1}+\cdots+m_{r}}e^{\left[m_{1}+\cdots+m+x\right]_{q}t}\nonumber \\
= & \sum_{n=0}^{\infty}E_{n,q}^{\left(r\right)}\left(x\right)\frac{t^{n}}{n!}.\nonumber
\end{align}

Thus, by (\ref{eq:10}), we get
\begin{align}
 & \int_{\Zp}\cdots\int_{\Zp}\left[x+y_{1}+\cdots+y_{r}\right]_{q}^{n}d\mu_{-1}\left(y\right)\label{eq:11}\\
= & E_{n,q}^{\left(r\right)}\left(x\right)\nonumber \\
= & \frac{1}{\left(1-q\right)^{n}}\sum_{l=1}^{n}\dbinom{n}{l}\left(-1\right)^{l}\left(\frac{2}{1+q^{l}}\right)^{r}\nonumber \\
= & 2^{r}\sum_{m_{1},\cdots,m_{r}=0}^{\infty}\left(-1\right)^{m_{1}+\cdots+m_{r}}\left[m_{1}+\cdots+m_{r}\right]_{q}^{n}.\nonumber
\end{align}

Let $w_{1},\, w_{2}\in\mathbb{N}$ with $w_{1}\equiv1$ (mod 2) and
$w_{2}\equiv1$ (mod 2). Then we observe that
\begin{align}
 & \int_{\Zp}\cdots\int_{\Zp}e^{\left[w_{1}\right]_{q}\left[w_{2}x+\frac{w_{2}}{w_{1}}\left(j_{1}+\cdots+j_{r}\right)+y_{1}+\cdots+y_{r}\right]_{q^{w_{1}}}t}d\mu_{-1}\left(y_{1}\right)\cdots d\mu_{-1}\left(y_{r}\right)\label{eq:12}\\
= & \int_{\Zp}\cdots\int_{\Zp}e^{\left[w_{1}w_{2}x+w_{2}\left(j_{1}+\cdots+j_{r}\right)+w_{1}\left(y_{1}+\cdots+y_{r}\right)\right]_{q}t}d\mu_{-1}\left(y_{1}\right)\cdots d\mu_{-1}\left(y_{r}\right)\nonumber \\
= & \lim_{N\rightarrow\infty}\sum_{y_{1},\cdots,y_{r}=0}^{p^{N}-1}e^{\left[w_{1}w_{2}x+w_{2}\left(j_{1}+\cdots+j_{r}\right)+w_{1}\left(y_{1}+\cdots+y_{r}\right)\right]_{q}t}\left(-1\right)^{y_{1}+\cdots+y_{r}}\nonumber \\
= & \lim_{N\rightarrow\infty}\sum_{y_{1},\cdots,y_{r}=0}^{w_{2}p^{N}-1}e^{\left[w_{1}w_{2}x+w_{2}\sum_{l=1}^{r}j_{l}+w_{1}\sum_{l=1}^{r}y_{l}\right]_{q}t}\left(-1\right)^{\sum_{l=1}^{r}y_{l}}\nonumber \\
= & \lim_{N\rightarrow\infty}\sum_{i_{1},\cdots,i_{r}=0}^{w_{2}-1}\sum_{y_{1},\cdots,y_{r}=0}^{p^{N}-1}e^{\left[w_{1}w_{2}x+w_{2}\sum_{l=1}^{r}j_{l}+w_{1}\sum_{l=1}^{r}(i_{l}+w_{2}y_{l})\right]_{q}t}\left(-1\right)^{\sum_{l=1}^{r}\left(i_{l}+w_{2}y_{l}\right)}.\nonumber
\end{align}

From (\ref{eq:12}), we have
\begin{align}
 & \sum_{j_{1},\cdots,j_{r}=0}^{w_{1}-1}\left(-1\right)^{j_{1}+\cdots+j_{r}}\label{eq:13}\\
 & \times\int_{\Zp}\cdots\int_{\Zp}e^{\left[w_{1}\right]_{q}\left[w_{2}x+\frac{w_{2}}{w_{1}}\left(j_{1}+\cdots+j_{r}\right)+y_{1}+\cdots+y_{r}\right]_{q^{w_{1}}}t}d\mu_{-1}\left(y_{1}\right)\cdots d\mu_{-1}\left(y_{r}\right)\nonumber \\
= & \lim_{N\rightarrow\infty}\left(\sum_{j_{1},\cdots,j_{r}=0}^{w_{1}-1}\sum_{i_{1},\cdots,i_{r}=0}^{w_{2}-1}\sum_{y_{1},\cdots,y_{r}=0}^{p^{N}-1}\left(-1\right)^{\sum_{l=1}^{r}\left(i_{l}+j_{l}+y_{l}\right)}\right.\nonumber \\
 & \left.\times e^{\left[w_{1}w_{2}\left(x+\sum_{l=1}^{r}y_{l}\right)+w_{2}\sum_{l=1}^{r}j_{l}+w_{1}\sum_{l=1}^{r}i_{l}\right]_{q}t}\right).\nonumber
\end{align}

By the same method as (\ref{eq:13}), we get
\begin{align}
 & \sum_{j_{1},\cdots,j_{r}=0}^{w_{2}-1}\left(-1\right)^{j_{1}+\cdots+j_{r}}\label{eq:14}\\
 & \times\int_{\Zp}\cdots\int_{\Zp}e^{\left[w_{2}\right]_{q}\left[w_{1}x+\frac{w_{1}}{w_{2}}\left(j_{1}+\cdots+j_{r}\right)+y_{1}+\cdots+y_{r}\right]_{q^{w_{2}}}t}d\mu_{-1}\left(y_{1}\right)\cdots d\mu_{-1}\left(y_{r}\right)\nonumber \\
= & \lim_{N\rightarrow\infty}\left(\sum_{j_{1},\cdots,j_{r}=0}^{w_{2}-1}\sum_{i_{1},\cdots,i_{r}=0}^{w_{1}-1}\sum_{y_{1},\cdots,y_{r}=0}^{p^{N}-1}\left(-1\right)^{\sum_{l=1}^{r}\left(i_{l}+j_{l}+y_{l}\right)}\right.\nonumber \\
 & \left.\times e^{\left[w_{1}w_{2}\left(x+\sum_{l=1}^{r}y_{l}\right)+w_{1}\sum_{l=1}^{r}j_{l}+w_{2}\sum_{l=1}^{r}i_{l}\right]_{q}t}\right).\nonumber
\end{align}

Therefore, by (\ref{eq:13}) and (\ref{eq:14}), we obtain the following
theorem.
\begin{thm}
\label{thm:1} For $w_{1},\, w_{2}\in\mathbb{N}$ with $w_{1}\equiv1$
$\textnormal{(mod 2})$ and $w_{2}\equiv1$ $\textnormal{(mod 2})$,
we have
\begin{align*}
& \sum_{j_{1},\cdots,j_{r}=0}^{w_{1}-1}\left(-1\right)^{j_{1}+\cdots+j_{r}}\\
 & \times\int_{\Zp}\cdots\int_{\Zp}e^{\left[w_{1}\right]_{q}\left[w_{2}x+\frac{w_{2}}{w_{1}}\left(j_{1}+\cdots+j_{r}\right)+y_{1}+\cdots+y_{r}\right]_{q^{w_{1}}}t}d\mu_{-1}\left(y_{1}\right)\cdots d\mu_{-1}\left(y_{r}\right)\nonumber \\
= & \sum_{j_{1},\cdots,j_{r}=0}^{w_{2}-1}\left(-1\right)^{j_{1}+\cdots+j_{r}}\\
 & \times\int_{\Zp}\cdots\int_{\Zp}e^{\left[w_{2}\right]_{q}\left[w_{1}x+\frac{w_{1}}{w_{2}}\left(j_{1}+\cdots+j_{r}\right)+y_{1}+\cdots+y_{r}\right]_{q^{w_{2}}}t}d\mu_{-1}\left(y_{1}\right)\cdots d\mu_{-1}\left(y_{r}\right).
\end{align*}

\end{thm}
$\,$
\begin{cor}
\label{cor:2} For $n\ge0$, and $w_{1},\, w_{2}\in\mathbb{N}$ with
$w_{1}\equiv1$ $\textnormal{(mod 2})$ and $w_{2}\equiv1$ $\textnormal{(mod 2})$,
we have
\begin{align*}
 & \left[w_{1}\right]_{q}^{n}\sum_{j_{1},\cdots,j_{r}=0}^{w_{1}-1}\left(-1\right)^{\sum_{l=1}^{r}j_{l}}\left.\int_{\Zp}\cdots\int_{\Zp}\right\{ \\
 & \left.\left[w_{2}x+\frac{w_{2}}{w_{1}}\left(j_{1}+\cdots+j_{r}\right)+\left(y_{1}+\cdots+y_{r}\right)\right]_{q^{w_{1}}}^{n}\right\} d\mu_{-1}\left(y_{1}\right)\cdots d\mu_{-1}\left(y_{r}\right)\\
= & \left[w_{2}\right]_{q}^{n}\sum_{j_{1},\cdots,j_{r}=0}^{w_{2}-1}\left(-1\right)^{\sum_{l=1}^{r}j_{l}}\left.\int_{\Zp}\cdots\int_{\Zp}\right\{ \\
 & \left.\left[w_{1}x+\frac{w_{1}}{w_{2}}\left(j_{1}+\cdots+j_{r}\right)+\left(y_{1}+\cdots+y_{r}\right)\right]_{q^{w_{2}}}^{n}\right\} d\mu_{-1}\left(y_{1}\right)\cdots d\mu_{-1}\left(y_{r}\right).
\end{align*}

\end{cor}
$\,$

Therefore, by (\ref{eq:11}) and Corollary \ref{cor:2}, we obtain
the following theorem.
\begin{thm}
\label{thm:3} For $n\ge0$, and $w_{1},\, w_{2}\in\mathbb{N}$ with
$w_{1}\equiv1$ $\textnormal{(mod 2})$ and $w_{2}\equiv1$ $\textnormal{(mod 2})$,
we have
\begin{align*}
 & \left[w_{1}\right]_{q}^{n}\sum_{j_{1},\cdots,j_{r}=0}^{w_{1}-1}\left(-1\right)^{\sum_{l=1}^{r}j_{l}}E_{n,q^{w_{1}}}^{\left(r\right)}\left(w_{2}x+\frac{w_{2}}{w_{1}}\left(j_{1}+\cdots+j_{r}\right)\right)\\
= & \left[w_{2}\right]_{q}^{n}\sum_{j_{1},\cdots,j_{r}=0}^{w_{2}-1}\left(-1\right)^{\sum_{l=1}^{r}j_{l}}E_{n,q^{w_{2}}}^{\left(r\right)}\left(w_{1}x+\frac{w_{1}}{w_{2}}\left(j_{1}+\cdots+j_{r}\right)\right).
\end{align*}

\end{thm}
$\,$

From (\ref{eq:11}), we can derive the following equation (\ref{eq:15})
:

\begin{align}
 & \int_{\Zp}\cdots\int_{\Zp}\left[w_{2}x+\frac{w_{2}}{w_{1}}\left(j_{1}+\cdots+j_{r}\right)+\left(y_{1}+\cdots+y_{r}\right)\right]_{q^{w_{1}}}^{n}d\mu_{-1}\left(y_{1}\right)\cdots d\mu_{-1}\left(y_{r}\right)\label{eq:15}\\
= & \sum_{i=0}^{n}\dbinom{n}{i}\left(\frac{\left[w_{2}\right]_{q}}{\left[w_{1}\right]_{q}}\right)^{i}\left[j_{1}+\cdots+j_{r}\right]_{q^{w_{2}}}^{i}q^{w_{2}\left(n-i\right)\sum_{l=1}^{r}j_{l}}\nonumber \\
 & \times\int_{\Zp}\cdots\int_{\Zp}\left[w_{2}x+\sum_{l=1}^{r}y_{l}\right]_{q}^{n-i}d\mu_{-1}\left(y_{1}\right)\cdots d\mu_{-1}\left(y_{r}\right)\nonumber \\
= & \sum_{i=0}^{n}\dbinom{n}{i}\left(\frac{\left[w_{2}\right]_{q}}{\left[w_{1}\right]_{q}}\right)^{i}\left[j_{1}+\cdots+j_{r}\right]_{q^{w_{2}}}^{i}q^{w_{2}\left(n-i\right)\sum_{l=1}^{r}j_{l}}E_{n-i,q^{w_{1}}}^{\left(r\right)}\left(w_{2}x\right).\nonumber
\end{align}

By (\ref{eq:15}), we get
\begin{align}
 & \left[w_{1}\right]_{q}^{n}\sum_{j_{1},\cdots,j_{r}=0}^{w_{1}-1}\left(-1\right)^{\sum_{l=1}^{r}j_{l}}\label{eq:16}\\
 & \times\int_{\Zp}\cdots\int_{\Zp}\left[w_{2}x+\frac{w_{2}}{w_{1}}\sum_{l=1}^{r}j_{l}+\sum_{l=1}^{r}y_{l}\right]_{q^{w_{1}}}^{n}d\mu_{-1}\left(y_{1}\right)\cdots d\mu_{-1}\left(y_{r}\right)\nonumber \\
= & \sum_{j_{1},\cdots,j_{r}=0}^{w_{1}-1}\left(-1\right)^{\sum_{l=1}^{r}j_{l}}\sum_{i=0}^{n}\dbinom{n}{i}\left[w_{2}\right]_{q}^{i}\left[w_{1}\right]_{q}^{n-i}\nonumber \\
 & \times\left[j_{1}+\cdots+j_{r}\right]_{q^{w_{2}}}^{i}q^{w_{2}\left(n-i\right)\sum_{l=1}^{r}j_{l}}E_{n-i,q^{w_{1}}}^{\left(r\right)}\left(w_{2}x\right)\nonumber \\
= & \sum_{i=0}^{n}\dbinom{n}{i}\left[w_{1}\right]_{q}^{n-i}\left[w_{2}\right]_{q}^{i}E_{n-i,q^{w_{1}}}^{\left(r\right)}\left(w_{2}x\right)\nonumber \\
 & \times\sum_{j_{1},\cdots,j_{r}=0}^{w_{1}-1}\left(-1\right)^{j_{1}+\cdots+j_{r}}\left[j_{1}+\cdots+j_{r}\right]_{q^{w_{2}}}^{i}q^{w_{2}\left(n-i\right)\sum_{l=1}^{r}j_{l}}\nonumber \\
= & \sum_{i=0}^{n}\dbinom{n}{i}\left[w_{1}\right]_{q}^{n-i}\left[w_{2}\right]_{q}^{i}T_{n,i,q^{w_{2}}}^{\left(r\right)}\left(w_{1}\right)E_{n-i,q^{w_{1}}}^{\left(r\right)}\left(w_{2}x\right),\nonumber
\end{align}
where
\begin{equation}
T_{n,i,q}^{\left(r\right)}\left(w\right)=\sum_{j_{1},\cdots,j_{r}=0}^{w-1}\left(-1\right)^{j_{1}+\cdots+j_{r}}q^{\left(n-i\right)\left(j_{1}+\cdots+j_{r}\right)}\left[j_{1}+\cdots+j_{r}\right]_{q}^{i}.\label{eq:17}
\end{equation}

By the same method as (\ref{eq:16}), we see that
\begin{align}
 & \left[w_{2}\right]_{q}^{n}\sum_{j_{1},\cdots,j_{r}=0}^{w_{2}-1}\left(-1\right)^{j_{1}+\cdots+j_{r}}\label{eq:18}\\
 & \times\int_{\Zp}\cdots\int_{\Zp}\left[w_{1}x+\frac{w_{1}}{w_{2}}\left(j_{1}+\cdots+j_{r}\right)+\left(y_{1}+\cdots+y_{r}\right)\right]_{q^{w_{2}}}^{n}d\mu_{-1}\left(y_{1}\right)\cdots d\mu_{-1}\left(y_{r}\right)\nonumber \\
= & \sum_{i=0}^{n}\dbinom{n}{i}\left[w_{1}\right]_{q}^{i}\left[w_{2}\right]_{q}^{n-i}T_{n,i,q^{w_{1}}}^{\left(r\right)}\left(w_{2}\right)E_{n-i,q^{w_{2}}}^{\left(r\right)}\left(w_{1}x\right).\nonumber
\end{align}

Therefore, by Corollary (\ref{cor:2}), (\ref{eq:16}), (\ref{eq:17})
and (\ref{eq:18}), we obtain the following theorem.
\begin{thm}
\label{thm:4} For $n\ge0$, and $w_{1},\, w_{2}\in\mathbb{N}$ with
$w_{1}\equiv1$ $\textnormal{(mod 2})$ and $w_{2}\equiv1$ $\textnormal{(mod 2})$,
we have
\begin{align*}
 & \sum_{i=0}^{n}\dbinom{n}{i}\left[w_{1}\right]_{q}^{n-i}\left[w_{2}\right]_{q}^{i}T_{n,i,q^{w_{2}}}^{\left(r\right)}\left(w_{1}\right)E_{n-i,q^{w_{1}}}^{\left(r\right)}\left(w_{2}x\right)\\
= & \sum_{i=0}^{n}\dbinom{n}{i}\left[w_{1}\right]_{q}^{i}\left[w_{2}\right]_{q}^{n-i}T_{n,i,q^{w_{1}}}^{\left(r\right)}\left(w_{2}\right)E_{n-i,q^{w_{2}}}^{\left(r\right)}\left(w_{1}x\right),
\end{align*}
where
\[
T_{n,i,q}^{\left(r\right)}\left(w\right)=\sum_{j_{1},\cdots,j_{r}=0}^{w-1}\left(-1\right)^{j_{1}+\cdots+j_{r}}q^{\left(n-i\right)\left(j_{1}+\cdots+j_{r}\right)}\left[j_{1}+\cdots+j_{r}\right]_{q}^{i}.
\]

\end{thm}


\bigskip
ACKNOWLEDGEMENTS. This work was supported by the National Research Foundation of Korea(NRF) grant funded by the Korea government(MOE)\\
(No.2012R1A1A2003786 ).
\bigskip


\bibliographystyle{amsplain}
\nocite{*}
\bibliography{0101}

$\,$

\noun{Department of Mathematics, Sogang University, Seoul 121-742,
Republic of Korea}

\emph{E-mail address : }\texttt{dskim@sogang.ac.kr}

$\,$

\noun{Department of Mathematics, Kwangwoon University, Seoul 139-701,
Republic of Korea}

\emph{E-mail address :} \texttt{tkkim@kw.ac.kr}
\end{document}